\documentclass[12pt,english]{article}
\usepackage[T1]{fontenc}
\usepackage[latin9]{inputenc}
\usepackage[a4paper]{geometry}
\usepackage{fancyhdr}
\usepackage{verbatim}
\usepackage{amsbsy}
\usepackage{graphicx}
\usepackage{setspace}
\usepackage{amssymb}
\usepackage{esint}
\makeatletter

\newcommand{\ds}{\displaystyle }

\newtheorem{mytheorem}{Theorem}
\newtheorem{mylemma}{Lemma}

\newcommand{\myQED}{\null\hfill $\square$}

\date{\mbox{\small  August 6, 2014}}

\makeatother

\usepackage{babel}

\begin{document}
\global\long\def\myargmax{\mathop{\mbox{arg\ensuremath{\,}max}}}

\global\long\def\HoSysSave{{\displaystyle \left.aH_{1}(t)\!\!\!\!\!\phantom{\int}\right|_{\abcd}}}
 \global\long\def\HoSys{{\displaystyle \left.\!\!\!\!\!\phantom{\int}\right|_{\abcd}}H_{1}(z)}
 \global\long\def\HoSysShort{\frac{d}{dt}{\displaystyle \left.\!\!\!\!\!\phantom{\int}\right|_{\abcd}\,}}
 \global\long\def\mysign{\mathop{\mbox{sign}}}
 \global\long\def\BRo{\mathcal{B}_{R_{0}}}
 \global\long\def\hR{\hat{R}}
 \global\long\def\hU{\hat{u}}
 \global\long\def\Atu#1{A(#1,t)}
 \global\long\def\Btu#1{B(#1,t)}

\title{Controlling of clock synchronization in WSNs: \\
structure of optimal solutions}

\author{Larisa Manita}

\maketitle
\vspace*{-4ex}

\noindent \begin{center}
\emph{National Research University Higher School of Economics}\\
\emph{ Moscow Institute of Electronics and Mathematics}\\
\emph{ Bolshoy Trehsviatitelskiy Per.~3, 109028 Moscow, Russia}
\par\end{center}

\noindent \begin{center}
\texttt{\vspace{1ex}
lmanita@hse.ru}
\par\end{center}
\begin{abstract}
\noindent Energy-saving optimization is very important for various
engineering problems related to modern distributed systems. We consider
here a control problem for a wireless sensor network with a single
time server node and a large number of client nodes. The problem is
to minimize a functional which accumulates clock synchronization errors
in the clients nodes and the energy consumption of the server over
some time interval $[0,T]$. The control function $u=u(t)$, $0\leq u(t)\leq u_{1}$,
corresponds to the power of the server node transmitting synchronization
signals to the clients. For all possible parameter values we find
the structure of optimal trajectories. We show that for sufficiently
large $u_{1}$ the solutions contain singular arcs.

\smallskip 

\end{abstract}

\section{Model }

Power consumption, clock synchronization and optimization are very
popular topics in analysis of wireless sensor networks \cite{arxiv}--\cite{MAC-protokol}.
In the majority of modern papers their authors discuss and compare
communication protocols (see, for example, \cite{toulouse}), network
architectures (for example, \cite{poljaki}) and technical designs
by using numerical simulations or dynamical programming methods (e.g.,
\cite{MicroSoft}). In~the~pre\-sent talk we consider a mathematical
model related with large scale networks which nodes are equipped with
noisy non-perfect clocks \cite{anatoly}. The task of optimal clock
synchronization in such networks is reduced to the classical control
problem. Its functional is based on the trade-off between energy consumption
and mean-square synchronization error. This control problem demon\-stra\-tes
surprisingly deep connections with the theory of singular optimal
solutions \cite{urs}--\cite{Powers}. 

The network consists of a single server node (denoted by~$1$) and
$N$ client nodes (sensors) numbered as $2,\ldots,N+1$.

Let $x_{i}$ be a state of the node~$i$ having the meaning of a
local clock value at this node. The network evolves in time $t\in\mathbb{R}_{+}$
as follows.

1) The node $1$ is a time server with the perfect clock: \[
\frac{d\, x_{1}(t)}{dt}=v>0\]

2) The client nodes are equipped with non-perfect clocks with a random
Gaussian noise\[
\frac{d\, x_{j}(t)}{dt}=v+\sigma dW_{j}(t)+\mbox{synchronizing jumps},\quad\]
where $W_{j}(t)$, $j=2,\ldots,N+1$, are independent standard Wiener
processes, $\sigma>0$ corresponds to the strength of the noise and
{}``synchronizing jumps'' are explained below.

3) At random time moments the server node~$1$ sends messages to
randomly chosen client nodes, $u$ is the \emph{intensity} of the
Poissonian message flow issued from the server. The client $j$, $j=2,\ldots,N+1$,
that receives at time $\tau$ a message from the node~$1$ immediately
ajusts its clock to the current value of~$x_{1}$: \begin{eqnarray*}
x_{j}(\tau+0) & = & x_{1}(\tau),\\
x_{k}(\tau+0) & = & x_{k}(\tau),\quad\quad k\not=j.\end{eqnarray*}

Hence the client clocks $x_{j}(t)$, $t\geq0$, are stochastic processes
which interact with the time server. 

The function \[
R(t)=\mathsf{E}\frac{1}{N}\sum_{j=2}^{N+1}\left(x_{j}(t)-x_{1}(t)\right)^{2}\]
is a cumulative measure of desynchronization between the client and
server nodes. Here $\mathsf{E}$ stands for the expectation.

It was proved in \cite{anatoly,a-man-TS-arxiv} that the function
$R(t)$ satisfies the differential equation \[
\dot{R}=-uR+N\sigma^{2}\]

\section{Optimal control problem}

Consider the following optimal control problem\begin{equation}
\int_{0}^{T}\left(\alpha R(t)+\beta u(t)\right)dt\rightarrow\inf\label{integral}\end{equation}
\begin{equation}
\dot{R}\left(t\right)=-u(t)R(t)+N\sigma^{2}\label{eq:r-t}\end{equation}
\begin{equation}
R\left(0\right)=R_{0}\label{eq:R_0}\end{equation}
\begin{equation}
0\leq u(t)\leq u_{1}\label{eq:ogr u}\end{equation}
Here $\alpha,\,\beta$ are some positive constants. The control function
$u\left(t\right)$ corresponds to the power of the server node transmitting
synchronization signals to the clients. The functional (\ref{integral})
 accumulates clock synchronization errors in the clients nodes and
the energy consumption of the server over some time interval $[0,T]$. 

The admissible solutions to (\ref{integral})-(\ref{eq:ogr u}) are
absolutely continuous functions, the admissible controls belong to
$L^{\infty}\left[0,T\right]$. 

We prove that the problem (\ref{integral})-(\ref{eq:ogr u}) has
a unique solution. We find a structure of optimal control. We show
that optimal solutions may contain singular arcs.

\section{Existence of solution}

\begin{mylemma} \label{l-lemma-1}

For any $R_{0}$ and any parameter values $T,\,\alpha,\,\beta,\, N,\,\sigma^{2},\, u_{1}$
there exists a solution $(\hR(t),\hU(t))$ to the problem (\ref{integral})-(\ref{eq:ogr u}).

\end{mylemma}

\noindent \emph{Proof.} Let $\BRo$ denote the set of continuous functions
$R:\,[0,T]\rightarrow\mathbb{R}$ such that $R(0)=R_{0}$. Consider
the map $K:\, L^{\infty}\left[0,T\right]\rightarrow\BRo$ defined
as follows:\[
\left(Ku\right)\left(t\right)=R{}_{0}\exp\left(-\int_{0}^{t}u(\xi)d\xi\right)+N\sigma^{2}\int_{0}^{t}\exp\left(-\int_{s}^{t}u(\xi)d\xi\right)ds\]
 \begin{equation}
=:\,\,\Atu u+\Btu u.\label{eq:k}\end{equation}

This operator assigns to the control function $u$ the corresponding
solution $R$ of (\ref{integral})-(\ref{eq:ogr u}). 

\textbf{1.} ~ Let $\left\{ u^{(n)}(t)\right\} _{n=1}^{\infty}$ be
a minimizing sequence for the fuctional \[
\int_{0}^{T}\left(\alpha R(t)+\beta u(t)\right)\, dt,\]
i.e., \[
\int_{0}^{T}\left(\alpha Ku^{(n)}(t)+\beta u^{(n)}(t)\right)\, dt\rightarrow\inf_{u\in V}\left\{ \int_{0}^{T}\left(\alpha R(t)+\beta u(t)\right)\, dt\right\} ,\qquad(n\rightarrow\infty),\]
where $V=\left\{ v\in L^{\infty}\left[0,T\right]:\,\,0\leq v(t)\leq u_{1}\right\} $.
Recall that the space $L^{1}\left[0,T\right]$ is the adjoint space
to $L^{\infty}\left[0,T\right]$. By $\left\langle \phi,u\right\rangle $
we denote the value of the functional $\phi\in\left(L^{\infty}\left[0,T\right]\right)^{*}\cong L^{1}\left[0,T\right]$
at $u\in L^{\infty}\left[0,T\right]$: \[
\left\langle \phi,u\right\rangle =\int_{0}^{T}\phi(\xi)u(\xi)\, d\xi\,.\]
Since $u^{(n)}(t)\in\left[0,u_{1}\right]$, one can extract a weakly-$*$
converging in $L^{\infty}\left[0,T\right]$ subsequence $u^{(n_{k})}(t)$
by virtue of Banach-Alaoglu theorem. Without loss of generality one
can assume that $u^{(n)}$ weakly-$*$ converges to some $\hU\in L^{\infty}\left[0,T\right]$.
This means that for each $\rho\in L^{1}\left[0,T\right]$ one has
\begin{equation}
\int_{0}^{T}\rho(\xi)u^{(n)}(\xi)\, d\xi\rightarrow\int_{0}^{T}\rho(\xi)\hU(\xi)\, d\xi,\,\quad n\rightarrow\infty.\label{eq:star}\end{equation}

\textbf{2.} ~ Let us prove that the sequence $R^{(n)}(t):=Ku^{(n)}(t)$
converges pointwise to $\hR(t):=K\hU(t)$ as $n\rightarrow\infty$.

Further let $\phi_{s}^{t}(\xi):=-\mathbf{1}_{\left[s,t\right]}\left(\xi\right)=\left\{ \begin{array}{cc}
-1, & \xi\in[s,t],\\
0, & \xi\not\in[s,t].\end{array}\right.$ Taking $\rho(\xi)=\phi_{0}^{t}(\xi)$ in (\ref{eq:star}) we obtain
\[
\int_{0}^{t}u^{(n)}(\xi)\, d\xi\rightarrow\int_{0}^{t}\hU(\xi)\, d\xi,\quad\, n\rightarrow\infty,\]
hence \[
\Atu{u^{(n)}}\rightarrow\Atu{\hU},\quad n\rightarrow\infty\]
for each fixed $t$. Note that $\Btu{u^{(n)}}={\displaystyle N\sigma^{2}\int_{0}^{t}\exp\left\langle \phi_{s}^{t},u^{(n)}\right\rangle ds}$.
The functions $\exp\left\langle \phi_{s}^{t},u^{(n)}\right\rangle $
are uniformly bounded and pointwise convergent, hence Lebesgue's dominated
theorem yields the convergence \[
\Btu{u^{(n)}}\rightarrow\Btu{\hU},\quad n\rightarrow\infty\]
for each fixed $t$. So we established the required convergence.

\textbf{3.} ~ Let us show that $\hR(t)$ is a solution to (\ref{integral})--(\ref{eq:ogr u}).

Obviously $R^{(n)}(t)$ are uniformly bounded (this follows straightforward
from the explicit formula (\ref{eq:k})). Since they form a pointwise
convergent sequence, Lebesgue's dominated theorem yields \[
\int_{0}^{T}\alpha R^{(n)}(t)\, dt\rightarrow\int_{0}^{T}\alpha\hR(t)\, dt,\quad\, n\rightarrow\infty.\]
Moreover, due to weak-$*$ convergence, one has \[
\int_{0}^{T}\beta u^{(n)}(t)\, dt=\beta\int_{0}^{T}\phi_{0}^{T}(t)u^{(n)}(t)dt\rightarrow\beta\int_{0}^{T}\phi_{0}^{T}(t)\hU(t)dt=\beta\int_{0}^{T}\hU(t)dt.\]

This yields \[
\int_{0}^{T}\left(\alpha R{}^{(n)}(t)+\beta u^{(n)}(t)\right)\, dt\rightarrow\int_{0}^{T}\left(\alpha\hR(t)+\beta\hU(t)\right)\, dt.\]

Thus $(\hR(t),\hU(t))$ is an optimal solution to (\ref{integral})--(\ref{eq:ogr u}).
\myQED

\section{Pontryagin maximum principle}

We will apply Pontryagin Maximum Principle \cite{Pontr} to the problem
(\ref{integral})-(\ref{eq:ogr u}). Let $\left(\widehat{R}\left(t\right),\widehat{u}\left(t\right)\right)$
be an optimal solution. Then there exist a constant $\lambda_{0}$
and a continuous function $\psi\left(t\right)$ such that for all
$t\in(0,T)$ we have \begin{equation}
H\left(\widehat{R}\left(t\right),\psi\left(t\right),\widehat{u}\left(t\right)\right)=\max_{0\leq u\leq u_{1}}H\left(\widehat{R}\left(t\right),\psi\left(t\right),u\right)\label{eq:usl-max}\end{equation}
where the Hamiltonian function \[
H\left(R,\psi,u\right)=-\lambda_{0}\left(\alpha R+\beta u\right)+\psi\left(-uR+N\sigma^{2}\right)\]

Except at points of discontinuity of $\widehat{u}\left(t\right)$
\begin{equation}
\dot{\psi}\left(t\right)=-\frac{\partial H\left(\widehat{R}\left(t\right),\psi\left(t\right),\widehat{u}\left(t\right)\right)}{\partial R}=\lambda_{0}\alpha+\widehat{u}\left(t\right)\psi\label{eq:adjoint}\end{equation}
 And $\psi$ satisfies the following transversality condition \begin{equation}
\psi\left(T\right)=0\label{eq:transver_psi}\end{equation}
The function $\psi\left(t\right)$ is called an adjoint function.
The condition (\ref{eq:usl-max}) is called the \textbf{maximum condition.} 

The dynamics equation (\ref{eq:r-t}) and the adjoint equation (\ref{eq:adjoint})
form a Hamiltonian system 

\begin{equation}
\begin{array}{rclrcl}
\dot{\psi} & = & \lambda_{0}\alpha+\widehat{u}\left(t\right)\psi\\
\dot{R} & = & -\widehat{u}\left(t\right)R+N\sigma^{2}\end{array}\label{ham-sys-1-1-1}\end{equation}
where $\widehat{u}\left(t\right)$ satisfies the maximum condition
(\ref{eq:usl-max}). The solutions $\left(R\left(t\right),\psi\left(t\right)\right)$
of (\ref{ham-sys-1-1-1}) are called extremals. If $\lambda_{0}\neq0$,
we say that $\left(R\left(t\right),\psi\left(t\right)\right)$ is
normal. One can show \cite{poljaki} that in the problem (\ref{integral})-(\ref{eq:ogr u})
every extremal is normal. So we can put $\lambda_{0}=1$.

\section{Switching function and singular extremals}

Denote \begin{equation}
H_{0}\left(R,\psi\right)=-\alpha R+\psi N\sigma^{2},\quad H_{1}\left(R,\psi\right)=-\beta-R\psi\label{eq:H1}\end{equation}
 then $H=H_{0}+uH_{1}$. The Hamiltonian $H$ is linear in $u$. Hence
to maximize it over the interval $u\in\left[0,u_{1}\right]$ we need
to use boundary values depending on the sign of $H_{1}$. \begin{equation}
\hat{u}(t)=\left\{ \begin{array}{cc}
0, & \quad H_{1}\left(R(t),\psi(t)\right)<0\\
u_{1,} & \quad H_{1}\left(R(t),\psi(t)\right)>0\end{array}\right.\label{eq:u_hat}\end{equation}
The function $H_{1}$ is called a switching function.

Suppose that there exists an interval $\left(t_{1},t_{2}\right)$
such that \begin{equation}
H_{1}\left(R(t),\psi(t)\right)=0,\quad\forall t\in\left(t_{1},t_{2}\right)\label{eq:H_1=00003D0}\end{equation}
then the extremal $\left(R\left(t\right),\psi\left(t\right)\right),\quad t\in\left(t_{1},t_{2}\right),$
is called a \textbf{\emph{singular}} one. In this case we can't find
an optimal control from the maximum condition (\ref{eq:usl-max}).
We will differentiate the identity $H_{1}\left(R(t),\psi(t)\right)\equiv0$
by virtue of the Hamiltonian system~(\ref{ham-sys-1-1-1}) until
a control $u$ appears with a non-zero coefficient.

We say that a number $q$ is the order of the singular trajectory
iff \[
\left.\frac{\partial}{\partial u}\,\frac{d^{k}}{dt^{k}}\right|_{(\mbox{\ref{ham-sys-1-1-1}})}H_{1}(R,\psi)=0,\quad\quad k=0,\ldots,2q-1,\]
\[
\left.\frac{\partial}{\partial u}\,\frac{d^{2q}}{dt^{2q}}\right|_{(\mbox{\ref{ham-sys-1-1-1}})}H_{1}(R,\psi)\neq0\]
in some open neighborhood of the singular trajectory $(R(t),\psi(t))$. 

It is known that $q$ is an integer.

Singular solutions arise frequently in control problems \cite{urs}-\cite{asmda}
and are therefore of practical significance. We prove that for suffiently
large $u_{1}$ a singular control is realised in the problem (\ref{integral})-(\ref{eq:ogr u}). 

\begin{mylemma}\label{l-lemma-2}

Let \[
\sqrt{\frac{\alpha N\sigma^{2}}{\beta}}\leq u_{1}\]
 then in the problem (\ref{integral})-(\ref{eq:ogr u}) there exists
a singular extremal of order $1$ \begin{equation}
\hat{R_{s}}\left(t\right)\equiv\sqrt{N\sigma^{2}\frac{\beta}{\alpha}},\quad\psi_{s}\left(t\right)\equiv-\sqrt{\frac{\alpha\beta}{N\sigma^{2}}}\label{eq:sing_extr}\end{equation}
 and the corresponding singular control is \[
u_{s}=\sqrt{\frac{\alpha N\sigma^{2}}{\beta}}\]

\end{mylemma}

\noindent \emph{Proof. }Assume that (\ref{eq:H_1=00003D0}) holds.
We will differentiate this identity along the extremal with respect
to $t$: \begin{eqnarray}
\left.\frac{d}{dt}\right|_{(\mbox{\ref{ham-sys-1-1-1}})}H_{1}(R\left(t\right),\psi\left(t\right))=0\quad & \Rightarrow & \quad-N\sigma^{2}\psi\left(t\right)-\alpha R\left(t\right)=0\label{eq:d-dt_H}\\
\left.\frac{d^{2}}{dt^{2}}\right|_{(\mbox{\ref{ham-sys-1-1-1}})}H_{1}(R\left(t\right),\psi\left(t\right))=0\quad & \Rightarrow & \quad u\left(\alpha R\left(t\right)-N\sigma^{2}\psi\left(t\right)\right)-2\alpha N\sigma^{2}=0\quad\,\label{eq:eq:d-dt-2_H}\end{eqnarray}

\noindent From (\ref{eq:H_1=00003D0}) and (\ref{eq:d-dt_H}) we have
\begin{equation}
R\left(t\right)=\sqrt{N\sigma^{2}\frac{\beta}{\alpha}},\quad\quad\psi\left(t\right)=-\sqrt{\frac{\alpha\beta}{N\sigma^{2}}}\label{eq:ocob}\end{equation}
Substituting (\ref{eq:ocob}) in (\ref{eq:eq:d-dt-2_H}) we obtain
\[
2\sqrt{N\sigma^{2}\alpha\beta}\cdot u-2\alpha N\sigma^{2}=0\]
Thus \[
R\left(t\right)\equiv\sqrt{N\sigma^{2}\frac{\beta}{\alpha}},\quad\quad\psi\left(t\right)\equiv-\sqrt{\frac{\alpha\beta}{N\sigma^{2}}}\]
is a singular extremal of order 1 and $u_{s}=\sqrt{\frac{\alpha N\sigma^{2}}{\beta}}$
is the corresponding singular control.

Note that if $\sqrt{\frac{\alpha N\sigma^{2}}{\beta}}>u_{1}$ then
$u_{s}$ does not satisfy the condition $0\leq u(t)\leq u_{1}$ hence
optimal solutions to the problem (\ref{integral})-(\ref{eq:ogr u})
are nonsingular. \myQED 

Recall the well-known generalized Legendre-Clebsch condition \cite{urs},
the necessary condition for optimality of the singular extremal of
order $1$: \[
\frac{\partial}{\partial u}\,\frac{d^{2}}{dt^{2}}H_{1}(\widehat{R}\left(t\right),\psi\left(t\right))\geq0\]

We see that this condition holds in our problem. One can show that
any concatenation of the singular control with a bang control $u=0$
or $u=u_{1}$ satisfies the necessary conditions of optimality \cite{urs}.

From the transversality condition (\ref{eq:transver_psi}) it is easily
seen that on the final time interval the optimal control $\widehat{u}\left(t\right)$
in the problem (\ref{integral})-(\ref{eq:ogr u})\emph{ }is nonsingular.
Namely, for all initial condition $R_{0}$ and for all parameter values
$\alpha,\,\beta,\, N,\,\sigma^{2},\, u_{1}$ we have the following
result.

\begin{mylemma} \label{l-lemma-3}

There exists $\varepsilon>0$ such that $\widehat{u}\left(t\right)=0$
for all $t\in\left(T-\varepsilon,T\right)$ . 

\end{mylemma}

\noindent \emph{Proof. }Using the transversality condition (\ref{eq:transver_psi})
we obtain $H_{1}(\widehat{R}\left(T\right),\psi\left(T\right))=-\beta<0$.
The continuity of the switching function $H_{1}$ implies that \[
H_{1}(\widehat{R}\left(t\right),\psi\left(t\right))<0\quad\forall t\in\left(T-\varepsilon,T\right)\]
 for some $\varepsilon>0$. The maximum condition (\ref{eq:usl-max})
yields $\widehat{u}\left(t\right)=0,\,\,\, t\in\left(T-\varepsilon,T\right)$.
\quad \null  \quad  \myQED

\section{The orbits of the Pontryagin maximum principle system }

Consider the behaviour of the extremals on the plane $\left(R,\psi\right)$.
Let $\Gamma$ be a switching curve, that is, a set of point such that
$H_{1}\left(R,\psi\right)=0$. By~(\ref{eq:H1}) we have $\Gamma=\left\{ \left.\left(R,\psi\right)\right|\beta+R\psi=0\right\} $.
We are interested in the domain $\left\{ \left(R,\psi\right):\,\, R>0\right\} $.
Denote \[
\Gamma^{+}=\Gamma\cap\left\{ \left(R,\psi\right):\,\, R>0\right\} \]
Above $\Gamma^{+}$ the optimal control $\hat{u}$ equals $0$, below
$\Gamma^{+}$ the optimal control $\hat{u}$ equals~$u_{1}$ (see~(\ref{eq:u_hat})).
Let $u=0$ then the Hamiltonian system (\ref{ham-sys-1-1-1}) has
the form \begin{equation}
\dot{R}=N\sigma^{2},\quad\dot{\psi}=\alpha\label{eq:PS_0}\end{equation}
The general solution of (\ref{eq:PS_0}) is \[
R\left(t\right)=N\sigma^{2}t+C_{1},\quad\psi\left(t\right)=\alpha t+C_{2}\]
 On the plane $\left(R,\psi\right)$ the orbits of the system (\ref{eq:PS_0})
are straight lines \[
\psi=\frac{\alpha}{N\sigma^{2}}R+C_{3}\]

Let $u=u_{1}$ than the Hamiltonian system (\ref{ham-sys-1-1-1})
has the form \begin{equation}
\dot{R}=-u_{1}R+N\sigma^{2},\quad\dot{\psi}=\alpha+u_{1}\psi\label{eq:PS_u1}\end{equation}
The general solution of (\ref{eq:PS_u1}) is \[
R(t)=\widetilde{C}e^{-u_{1}t}+\frac{N\sigma^{2}}{u_{1}},\quad\psi\left(t\right)=\widetilde{w}e^{u_{1}t}-\frac{\alpha}{u_{1}}\]

On the plane $\left(R,\psi\right)$ if $\widetilde{C}\neq0,\,\,\widetilde{w}\neq0$,
the orbits of the system (\ref{eq:PS_u1}) are hyperbolas \[
\left|\alpha+\psi u_{1}\right|\cdot\left|N\sigma^{2}-u_{1}R\right|=\omega\]
If $\widetilde{C}=0,\,\,\widetilde{w}\neq0$, the orbit is the straight
line $R=\frac{N\sigma^{2}}{u_{1}}$, directed upward if $\widetilde{w}>0$
or downward if $\widetilde{w}<0$. If $\widetilde{w}=0$, the orbit
is the straight line $\psi=-\frac{\alpha}{u_{1}}$, directed to the
left if $\widetilde{C}>0$ or to the right if $\widetilde{C}<0$.
If $\widetilde{C}=0,\,\,\widetilde{w}=0$ , the point $\left(\frac{N\sigma^{2}}{u_{1}},-\frac{\alpha}{u_{1}}\right)$
is the stationary orbit. 

\vspace*{1.5cm}


\parbox{1.\textwidth}{ 

\null\vspace*{-3cm}

\begin{center}
\null\hspace*{-1cm}\includegraphics[width=1.1\textwidth]{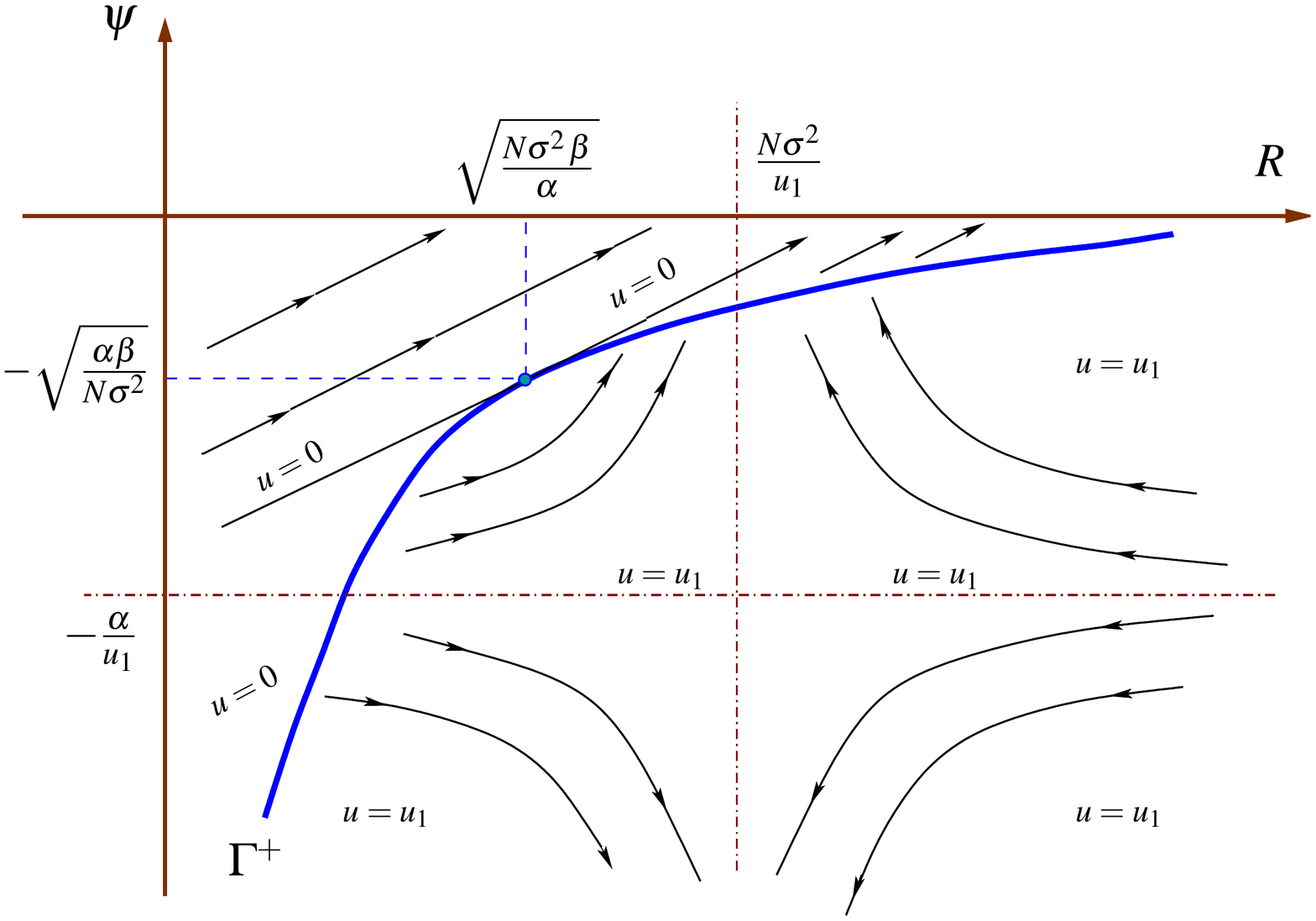}
\par\end{center}

\vspace*{-1.5cm}

\begin{center}
\textbf{Fig 1.} Orbits in the nonsingular case: ${\displaystyle \sqrt{\frac{\alpha N\sigma^{2}}{\beta}}}\,>u_{1}$ 
\par\end{center}

\vspace*{0.5cm}}  

\parbox{1.\textwidth}{ 

\null\vspace*{-1.5cm}

\begin{center}
\null\hspace*{-1cm}\includegraphics[width=1.1\textwidth]{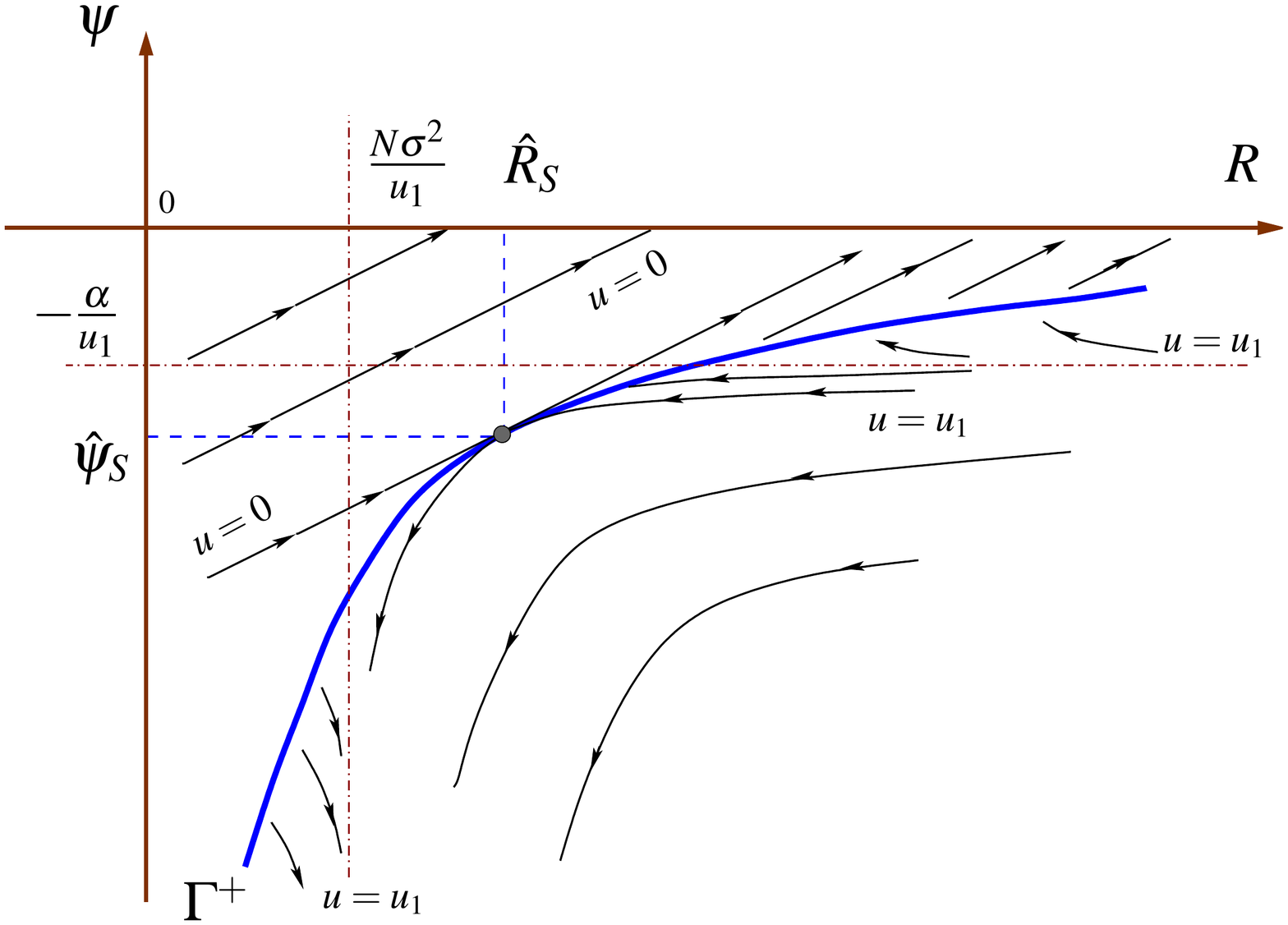}
\par\end{center}

\vspace*{-1.15cm}

\begin{center}
\textbf{Fig 2.} Orbits in the singular case: ${\displaystyle \sqrt{\frac{\alpha N\sigma^{2}}{\beta}}}\,\leq u_{1}$ 
\par\end{center}

\vspace*{0.03cm}}  

\noindent \textbf{Remark.} On these figures we don't show trajectories
$(R(t),\psi(t))$ with  $\psi\left(0\right)>0$ because they cannot
satisfy the transversality condition. 

\vspace*{0.5cm}


Note that in the case ${\displaystyle \sqrt{\alpha N\sigma^{2}/\beta\,}}\,\leq u_{1}$
two extremals go out of the singular point $\left(\sqrt{N\sigma^{2}\frac{\beta}{\alpha}},-\sqrt{\frac{\alpha\beta}{N\sigma^{2}}}\right)$
(with $u=0$ and $u=u_{1}$). But only one extremal (going of the
singular point) satisfies the transversality condition~(\ref{eq:transver_psi}). 

Thus for any $R_{0}\geq0$ there exists a unique extremal such that
$R\left(0\right)=R_{0},$ $\psi\left(T\right)=0$. Since we prove
that a solution to problem (\ref{integral})-(\ref{eq:ogr u}) exists
hence the constructed extremals are optimal. 

To summarize the above analysis in the next two sections we consider
separately the nonsingular and singular cases. In each case we provide
a plot with optimal solutions and state a conclusion on the structure
of the optimal control $\hat{u}(t)$ (Theorems~\ref{my-Theorem-1}
and~\ref{my-Theorem-2}). It is interesting also to see how the structure
of $\hat{u}(t)$ depends on the parameter $R_{0}$ and~$T$. The
answer is presented on Figures~4 and~6.

\newpage

\section{Optimal solutions. Nonsingular case}

\parbox{1.\textwidth}{ 

\null\vspace*{-2cm}

\begin{center}
\null\hspace*{-1.5 cm} \includegraphics[width=1.18\textwidth]{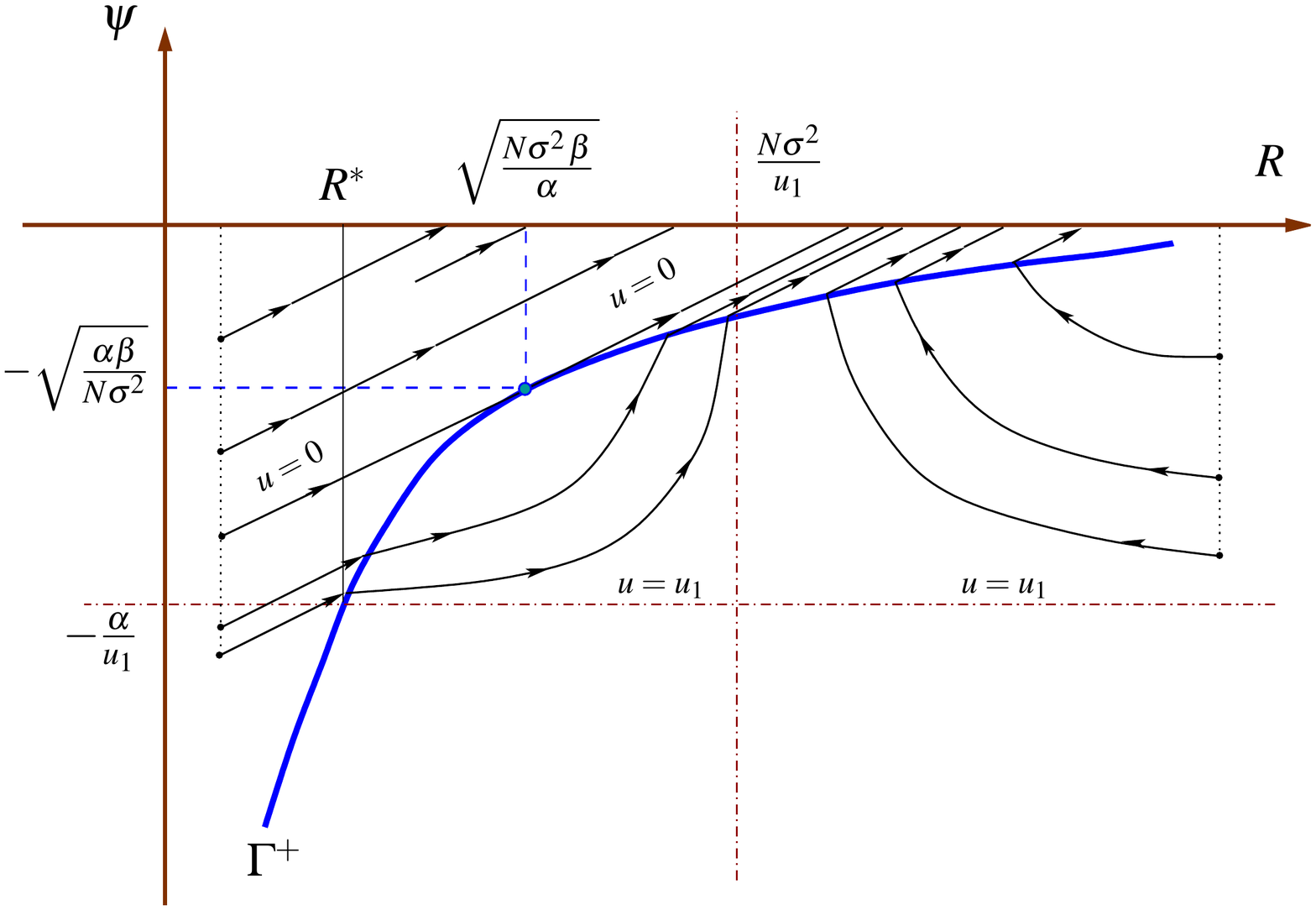}
\par\end{center}

\vspace*{-2cm}

\begin{center}
\textbf{Fig 3.} Optimal solutions for different values of the problem
parameters. \\
\emph{Nonsingular case}. 
\par\end{center}

\vspace*{0.5cm}}  

\begin{mytheorem}\label{my-Theorem-1}

Let $\ds\sqrt{\frac{\alpha N\sigma^{2}}{\beta}}>u_{1}$, that is,
optimal solutions are nonsingular \\
(Lemma 2). Then, depending of values $R\left(0\right)$ and $T$,
the optimal control $\hat{u}(t)$ has one of the following forms \[
1.1.\quad\,\,\hat{u}(t)=0,\,\, t\in(0,T)\]
\[
1.2.\quad\,\,\hat{u}(t)=\left\{ \begin{array}{cc}
u_{1}, & t\in(0,t_{1})\\
0, & t\in(t_{1},T)\end{array}\right.\]
 \[
1.3.\quad\,\,\hat{u}(t)=\left\{ \begin{array}{cc}
0, & t\in(0,t_{1})\\
u_{1}, & t\in(t_{1},t_{2})\\
0, & t\in(t_{2},T)\end{array}\right.\]
 i.e., the optimal control switches between $u=0$ and $u=u_{1}$
and the number of switchings does not exceed 2. 

\end{mytheorem} 

\newpage 

The Fig.~4 shows how the \emph{structure} of optimal controls $\hat{u}=\hat{u}(t)$,
$t\in[0,T]$, depends on $T$ and on the initial value~$R(0)$.

\parbox{1.\textwidth}{ 

\null\vspace*{-.2cm}

\hspace*{-1cm}\includegraphics[width=1.15\textwidth]{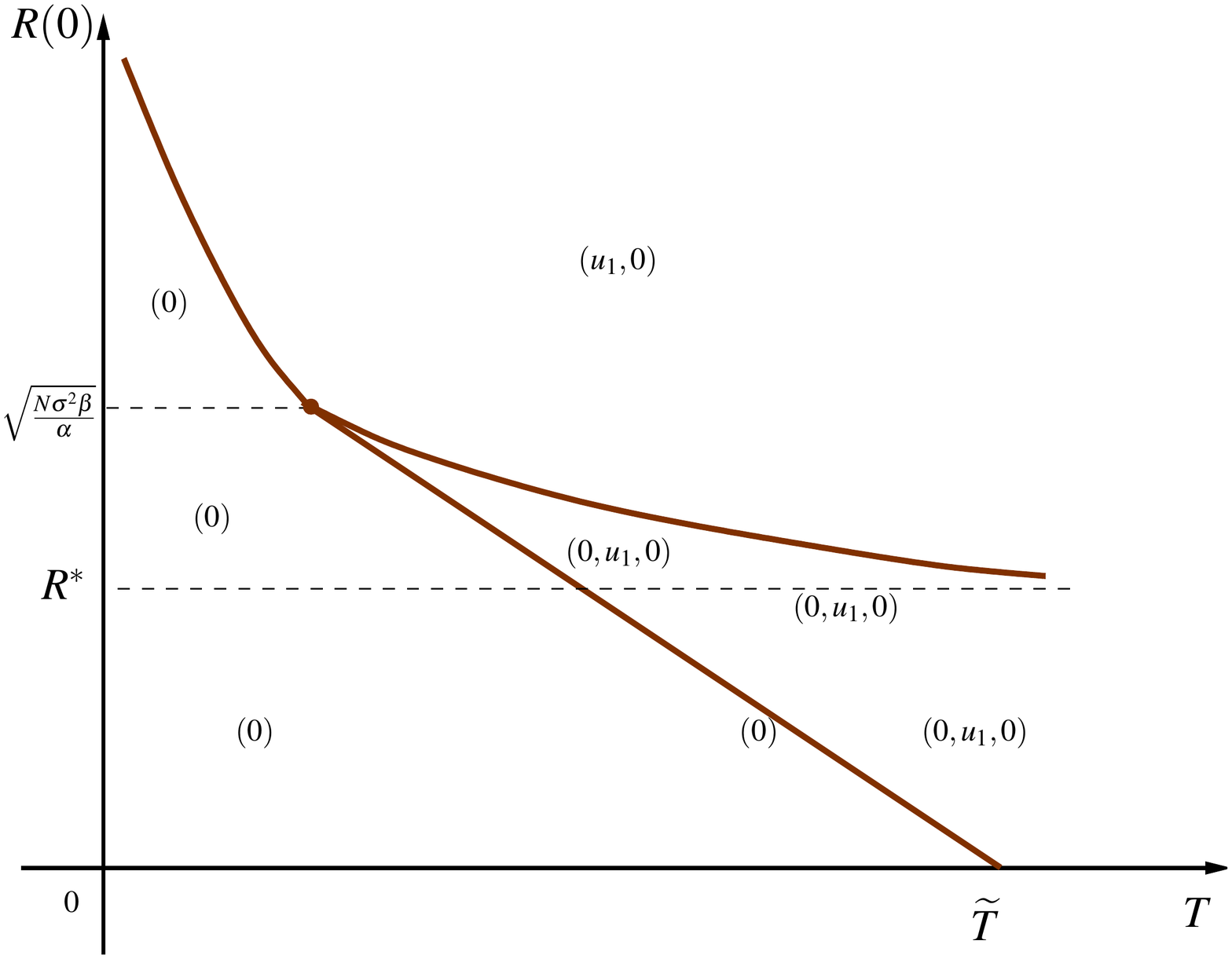}

\vspace*{-1cm}

\begin{center}
\textbf{Fig 4.} ~ 
\par\end{center}

\vspace*{0.05cm}}  

\noindent Let $(\theta,\rho)$ be some point on the plane $\left(T,R(0)\right)$.
Assume that $(\theta,\rho)$ belongs to a domain labeled, \emph{for
example}, by $(a,b,c)$. This means that for the optimal control problem
with $T=\theta$ and $R(0)=\rho$ the optimal control function $\hat{u}(t)$
has the following form \[
\hat{u}(t)=\left\{ \begin{array}{rcl}
a, &  & t\in(0,\tau_{1}),\\
b, &  & t\in(\tau_{1},\tau_{2}),\\
c, &  & t\in(\tau_{2},\theta).\end{array}\right.\]
Here $\tau_{1}$ and $\tau_{2}$ are some numbers satisfying the condition
$0<\tau_{1}<\tau_{2}<\theta$. The numbers $\tau_{1}$ and $\tau_{2}$
depend on $(\theta,\rho)$ and on all parameters $(\alpha,\beta,N,\sigma)$
of the model. For points $(\theta,\rho)$ in the domain labeled by
$(0)$ we have $\hat{u}(t)=0$ for all $t\in[0,T]$.

\section{Optimal Solutions. Singular case}

\parbox{1.\textwidth}{ 

\null\vspace*{-.2cm}

\begin{center}
\includegraphics[width=0.98\textwidth]{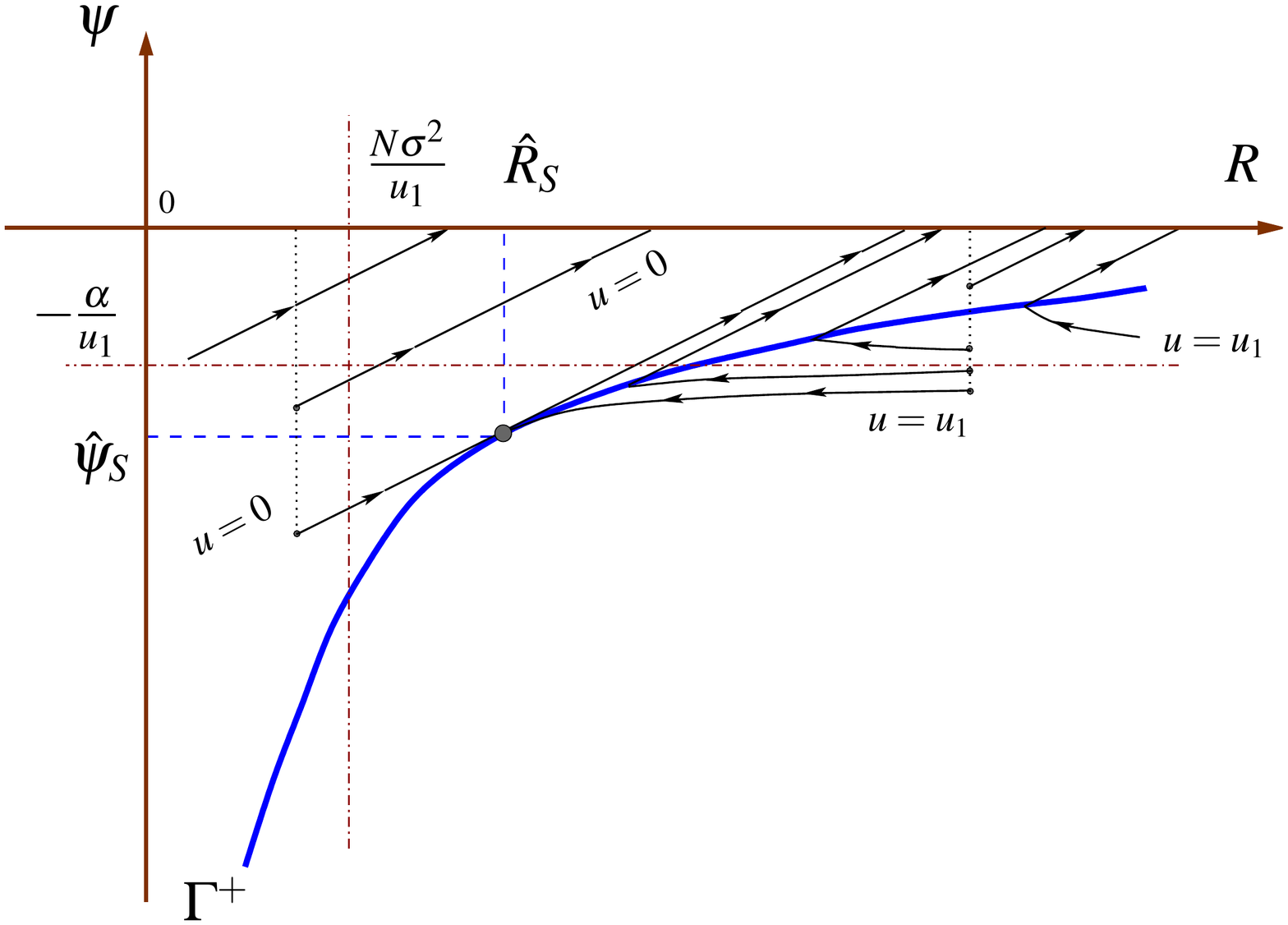}
\par\end{center}

\vspace*{-1cm}

\begin{center}
\textbf{Fig 5.} ~ Optimal solutions for different values of the model
parameters. \\
\emph{Singular case}. 
\par\end{center}

\vspace*{0.05cm}}  

\begin{mytheorem} \label{my-Theorem-2}

Let $\ds\sqrt{\frac{\alpha N\sigma^{2}}{\beta}}\leq u_{1}$. Then,
depending of values $R\left(0\right)$ and $T$, the optimal control
$\hat{u}(t)$ has one of the following forms \[
2.1.\quad\,\,\hat{u}(t)=0,\,\, t\in(0,T)\]
\[
2.2.\quad\,\,\hat{u}(t)=\left\{ \begin{array}{cc}
u_{1}, & t\in(0,t_{1})\\
0, & t\in(t_{1},T)\end{array}\right.\quad2.3.\quad\,\,\hat{u}(t)=\left\{ \begin{array}{cc}
u_{s}, & t\in(0,t_{1})\\
0, & t\in(t_{1},T)\end{array}\right.\]
 \[
2.4.\,\quad\,\hat{u}(t)=\left\{ \begin{array}{cc}
0, & t\in(0,t_{1})\\
u_{s}, & t\in(t_{1},t_{2})\\
0, & t\in(t_{2},T)\end{array}\right.,\quad2.5.\,\quad\,\hat{u}(t)=\left\{ \begin{array}{cc}
u_{1}, & t\in(0,t_{1})\\
u_{s}, & t\in(t_{1},t_{2})\\
0, & t\in(t_{2},T)\end{array}\right.\]
i.e., the number of control switchings does not exceed 2 and the optimal
solutions may contain the singular arcs (cases 2.3-2.5).

\end{mytheorem}

\parbox{1.\textwidth}{ 

\null\vspace*{-1.2cm}

\hspace*{-2cm}\includegraphics[width=1.15\textwidth]{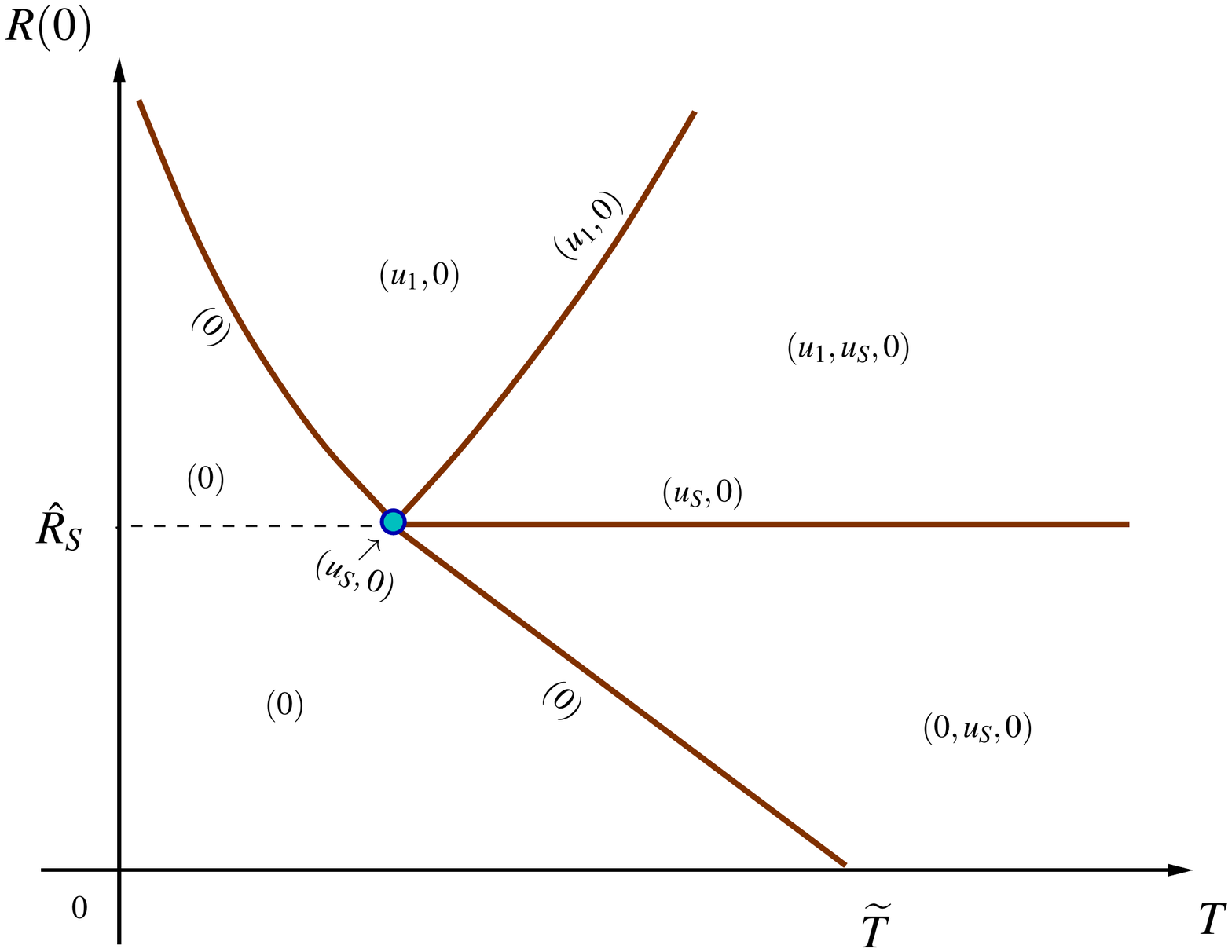}

\vspace*{-1cm}

\begin{center}
\textbf{Fig 6.} ~ 
\par\end{center}

\vspace*{0.05cm}}  

As it is seen from Fig.~6 in the \emph{singular case} on the plane
$\left(T,R(0)\right)$ we have more domains with different structures
of the optimal control $\hat{u}=\hat{u}(t)$. These additional domains
are labeled as $(u_{S},0)$ or $(a,u_{S},0)$. Note that on that intervals
$t\in\Delta$ where $\hat{u}(t)=u_{S}$ the function $\hat{R}(t)$
takes the constant value $\hat{R}_{S}$: \[
\hat{R}(t)=\hat{R}_{S},\qquad t\in\Delta.\]

\section{Conclusions}

We considered the control problem for wireless sensor networks with
a single time server node and a large number of client nodes. The
cost functional of this control problem accumulates clock synchronization
errors in the clients nodes and the energy consumption of the server
over some time interval $[0,T]$. For all possible parameter values
we found the structure of optimal control function. It was proved
that for any optimal solution $\widehat{R}\left(t\right)$ there exist
a time moment $\tau,$ $0\leq\tau<T$, such that $\hat{u}(t)=0,$
$t\in[\tau,T]$, i.e., the sending messages at times close to $T$
is not optimal. We showed that for sufficiently large $u_{1}$ the
optimal solutions contain singular arcs. We found conditions on the
model parameters under which different types of the optimal control
are realized.

We hope that our study of the energy-saving optimization will also
be usefull for analysis of other engineering problems related to modern
distributed systems. In future we plan to extend these results to
more general models.

\end{document}